\documentclass[11pt]{amsart}
\usepackage{amsfonts}
\usepackage{amsmath}
\usepackage{amsthm}
\usepackage{bm}
\usepackage{bbm}
\usepackage{graphicx}
\usepackage{float}
\usepackage{enumerate}
\usepackage{natbib}

\usepackage{pdfsync}

\usepackage{mathtools}

\usepackage{blindtext}
\usepackage[utf8]{inputenc}

\usepackage[hmarginratio=1:1,top=32mm,columnsep=20pt]{geometry}

\usepackage[T1]{fontenc} 
\usepackage{color}	
\usepackage[english]{babel} 

\usepackage{enumitem} 
\setlist[itemize]{noitemsep} 

\usepackage{leftidx}   

\usepackage{mathrsfs}

 \usepackage[all]{xy}

\begin{document}
\allowdisplaybreaks[4]
\newtheorem{theorem}{Theorem}
\newtheorem{lemma}{Lemma}
\newtheorem{pron}{Proposition}
\newtheorem{re}{Remark}
\newtheorem{thm}{Theorem}
\newtheorem{Corol}{Corollary}
\newtheorem{exam}{Example}
\newtheorem{defin}{Definition}
\newtheorem{remark}{Remark}
\newtheorem{property}{Property}
\newcommand{\blue}{\color{blue}}
\newcommand{\red}{\color{red}}
\newcommand{\la}{\frac{1}{\lambda}}
\newcommand{\sectemul}{\arabic{section}}
\renewcommand{\theequation}{\sectemul.\arabic{equation}}
\renewcommand{\thepron}{\sectemul.\arabic{pron}}
\renewcommand{\thelemma}{\sectemul.\arabic{lemma}}
\renewcommand{\there}{\sectemul.\arabic{re}}
\renewcommand{\thethm}{\sectemul.\arabic{thm}}
\renewcommand{\theCorol}{\sectemul.\arabic{Corol}}
\renewcommand{\theexam}{\sectemul.\arabic{exam}}
\renewcommand{\thedefin}{\sectemul.\arabic{defin}}
\renewcommand{\theremark}{\sectemul.\arabic{remark}}

\def\REF#1{\par\hangindent\parindent\indent\llap{#1\enspace}\ignorespaces}
\def\lo{\left}
\def\ro{\right}
\def\be{\begin{equation}}
\def\ee{\end{equation}}
\def\beq{\begin{eqnarray*}}
\def\eeq{\end{eqnarray*}}
\def\bea{\begin{eqnarray}}
\def\eea{\end{eqnarray}}
\def\o{\overline}
\newcommand{\bH}{\mathbf{H}}
\newcommand{\bg}{\mathbf{g}}
\newcommand{\ba}{\mathbf{a}}
\newcommand{\bb}{\mathbf{b}}
\newcommand{\bT}{\mathbf{T}}
\newcommand{\bh}{\mathbf{h}}
\newcommand{\bs}{\mathbf{s}}
\newcommand{\bY}{\mathbf{Y}}
\newcommand{\cid}{\stackrel{d}{\to}}
\newcommand{\cip}{\stackrel{P}{\to}}

\newcommand{\reals}{{\mathbb R}}
\newcommand{\bbr}{\reals}
\newcommand{\eid}{\overset{d}{=}}

\title{Cauchy, normal and correlations versus heavy tails}

\author{Hui Xu}
\address{Center for Applied Mathematics\\
Cornell University}
\email{hx223@cornell.edu}

\author{Joel E. Cohen}
\address{Laboratory of Populations, 
  The Rockefeller University \& Columbia University; Earth Institute \& Department of Statistics, Columbia University; 
  Department of Statistics, University of Chicago}
\email{cohen@rockefeller.edu}

\author{Richard Davis}
\address{Department of Statistics \\
  Columbia University}
\email{rdavis@stat.columbia.edu}

\author{Gennady Samorodnitsky$^*$}
\address{School of Operations Research and Information Engineering\\
Cornell University}
\email{gs18@cornell.edu}
  
\numberwithin{equation}{section}
 \thanks{$^*$The corresponding author.  This research was partially
   supported by NSF grants  DMS-2015379 (Davis) at Columbia and DMS-2015242 (Samorodnitsky) at Cornell.}

\subjclass{Primary  }
\keywords{heavy tails,  dependence, Cauchy, normal}

\begin{abstract}
A surprising result of \cite{pillai:meng:2016} showed that a 
transformation $\sum_{j=1}^n w_j X_j/Y_j$ of two iid 
centered normal random vectors, $(X_1,\ldots, X_n)$ and
$(Y_1,\ldots, Y_n)$, $n>1$, 
for any weights $0\leq w_j\leq 1$, $ j=1,\ldots, n$, $\sum_{j=1}^n w_j=1$,
has a Cauchy distribution regardless of any correlations within the normal vectors. 
The correlations appear to lose out in the competition with the heavy tails. 
To clarify how extensive this phenomenon is, we analyze two other transformations
of two iid centered normal random vectors.
These transformations are similar in spirit to the
transformation considered by \cite{pillai:meng:2016}. 
One transformation involves absolute values: 
$\sum_{j=1}^n w_j X_j/|Y_j|$.
The second involves randomly stopped Brownian motions:
$\sum_{j=1}^n w_j X_j\bigl(Y_j^{-2}\bigr)$, where
$\bigl\{\bigl( X_1(t),\ldots, X_n(t)\bigr), \, t\geq 0\bigr\},\ n>1,$ 
is a Brownian motion with positive variances;
$(Y_1,\ldots, Y_n)$ is a centered normal random
vector with the same law as $( X_1(1),\ldots, X_n(1))$ and independent of it; 
and $X(Y^{-2})$ is the value of the Brownian motion $X(t)$ evaluated
at the random time $t=Y^{-2}$. 
All three transformations result in a Cauchy distribution if the covariance
matrix of the normal components is diagonal, or if all the
correlations implied by the covariance matrix equal 
1. However, while the transformation \cite{pillai:meng:2016} considered 
produces a Cauchy distribution regardless of the normal covariance matrix.
the transformations we consider here
do not always produce a Cauchy distribution. 
The correlations between jointly normal random
variables are not always overwhelmed by the heaviness of the marginal tails. 
The mysteries of the connections between normal and Cauchy
laws remain to be understood. 
\end{abstract}

\maketitle

\section{Introduction} \label {sec:intro}
The normal distribution and the Cauchy distribution are among the most 
fundamental probability distributions with numerous applications in
many areas where uncertainty is present. 
They both arise as limits of suitably normalized sums of independent copies of random quantities. 
They both have nearly magical invariance properties. 
They contrast strongly in that the normal distribution has finite moments of all 
orders (even all exponential moments are finite), while the Cauchy distribution does
not even have a finite mean. It is well known that normal
and Cauchy random variables are
connected in simple ways. These connections sometimes appear
to have mysterious properties as well. 
This note aims to shed additional light to some of these mysteries.

By definition, the centered normal distribution with standard deviation 
$\sigma>0$ has the density
$$
\phi(x;\sigma) = \frac{1}{\sigma\sqrt{2\pi}} \exp\bigl\{
-x^2/(2\sigma^2)\bigr\}, \ x\in\bbr,
$$
while the Cauchy distribution with scale $\sigma>0$ has the density
\begin{equation} \label{e:dens.cauchy}
f(x;\sigma)=\frac{\sigma}{\pi(x^2+\sigma^2)}, \ x\in\bbr. 
\end{equation}
For both distributions, the case $\sigma=1$ is called standard. The
invariance properties of these distributions say that, if $X_1,\ldots,
X_n$ are iid copies of a centered normal $X$ and $\eid$ denotes equality
in distribution, then 
\begin{equation} \label{e:invariance}
\sum_{j=1}^n w_jX_j \eid X
\end{equation}
if $\sum_{j=1}^n w_j^2=1$, while if $X_1,\ldots,
X_n$ are iid copies of a Cauchy $X$, then \eqref{e:invariance} holds
if $\sum_{j=1}^n |w_j|=1$.

Versions of \eqref{e:invariance} may hold when the $X_1,\ldots, X_n$ 
are not independent. 
For example, recall that $X_1,\ldots, X_n$ are
jointly Cauchy if there is a finite symmetric measure $\Gamma$ (the
spectral measure) on the unit
sphere $S^n$ such that
\begin{equation} \label{e:mult.cauchy}
E\exp\biggl\{ i\sum_{j=1}^n \theta_jX_j\biggr\}
= \exp\biggl\{ -\int_{S^n} \biggl| \sum_{j=1}^d
    \theta_js_j\biggr|\,\Gamma(ds_1,\ldots, ds_n)\biggr\}.
\end{equation} 
Because $(X_1, \ldots,X_n)$ has a symmetric distribution,
the characteristic function is real-valued on the right side of \eqref{e:mult.cauchy}.
Then $X_1,\ldots, X_n$ have
identical marginals if $\int |s_j|\,\Gamma(d\bs)$ is the same for all
$j$; see \cite{samorodnitsky:taqqu:1994}. In particular, if the
spectral measure $\Gamma$ is concentrated on the positive and negative
quadrants of the unit sphere and 
$X_1,\ldots, X_n$ have the same marginals, then \eqref{e:invariance}
still holds as long as all $(w_j)$ are also of the same sign. 
In this scenario, $X_1,\ldots, X_n$ are generally dependent since
their independence requires $\Gamma$ to be concentrated at the points
where the unit sphere intersects the axes. 
In particular, if $\Gamma$ is concentrated at the points $\pm (1,\ldots,1)/n^{1/2}$,
then $X_1=\ldots=X_n$ a.s. 
This special case shows that a convex linear combination of iid Cauchy
random variables has the same law as a convex linear
combination of almost surely equal Cauchy random variables! 

\section{The Pillai and Meng result} \label{sec:pillai.meng}

An elementary calculation shows that, if
$X,Y$ are iid centered normal and $W$ is a standard Cauchy, then
\begin{equation} \label{e:normal.cauchy}
  X/Y\eid W.
\end{equation}
If $X$ and $Y$ are bivariate normal, then the probability density function of the ratio $X/Y$ is derived in \cite{Epham2006density}
and earlier sources they cite. 
The ratio has a Cauchy distribution if and only if $X$ and $Y$ are centered and uncorrelated.  For the remainder of this note, we will take $X$ and $Y$ to be centered and uncorrelated.

\cite{pillai:meng:2016} proved an amazing extension of \eqref{e:normal.cauchy}.
Let $(X_1,\ldots, X_n)$ and
$(Y_1,\ldots, Y_n),\ n>1,$ be iid multivariate centered normal random
vectors with positive variances and a covariance matrix $\Sigma$. Then
for any weights $0\leq w_j\leq 1$, $ j=1,\ldots, n, \,
\sum_{j=1}^n w_j=1$,
\begin{equation} \label{e:pillai.meng}
\sum_{j=1}^n w_j \frac{X_j}{Y_j}\eid W,
\end{equation}
where $W$ is a standard Cauchy.
Earlier, \cite{drton:xiao:2016} proved \eqref{e:pillai.meng}
when $n=2$. 

What is amazing in \eqref{e:pillai.meng}?
Each ratio $X_j/Y_j$ in the sum has the standard Cauchy
distribution by \eqref{e:normal.cauchy}. If the covariance matrix
$\Sigma$ is diagonal, then these Cauchy random variables are
independent, and \eqref{e:pillai.meng} is simply the invariance
property \eqref{e:invariance}. 
If, on the other hand, all the
correlations implied by the covariance matrix $\Sigma$ equal 1,
then all terms in \eqref{e:pillai.meng} are equal almost surely, and
\eqref{e:pillai.meng} is equivalent to the above-mentioned fact that 
a convex linear combination of equal Cauchy random variables has the same law as a
convex linear combination of iid Cauchy
random variables. The amazing part is that \eqref{e:pillai.meng} holds
for all other covariance matrices $\Sigma$. 

The only comparable property of non-independent Cauchy random variables 
that appears to be known at this point is the property of jointly
Cauchy random variables with the
spectral measure $\Gamma$ concentrated on the positive and negative
quadrants of the unit sphere. 
It is straightforward to show
that the terms in \eqref{e:pillai.meng} are jointly
Cauchy if and only if the normal random vectors consist of independent blocks
of equal normal random variables. 

Therefore,  
it is a special feature of the multivariate normal distribution
that 
the ratios of the components of iid centered normal random
vectors have the invariance property \eqref{e:invariance},
regardless of the correlations among the normal components.
\cite{pillai:meng:2016} speculate that ``the dependence among them'' 
(i.e., among the ratios of centered,  dependent normal random variables) 
``can be overwhelmed by
the heaviness of their marginal tails'' (of the ratios) ``in determining the
stochastic behavior of their linear combinations.'' This speculation
continues to generate interest and new results.
For example, \cite{cds:2020} give another example where heavy tails overwhelm
correlations among normal random variables.

Here we shall show, by two examples, that 
the correlations between jointly normal random
variables are not always overwhelmed by the heaviness of the marginal tails.
In the two subsequent sections, we analyze two other natural transformations,
in the spirit of \cite{pillai:meng:2016},
from the normal world to the Cauchy world. 

\section{Transformation 1: Absolute Values} \label{sec:absolute.val}

The symmetry of a centered normal random variable and
\eqref{e:normal.cauchy} imply that,
if $X,\ Y$ are iid centered normals and $W$ is a standard Cauchy, then
\begin{equation} \label{e:normal.cauchy.1}
X/|Y|\eid W.
\end{equation}
Following \cite{pillai:meng:2016}, let  $(X_1,\ldots, X_n)$ and
$(Y_1,\ldots, Y_n),\ n>1$, be iid multivariate centered normal random
vectors with positive variances 
and a covariance matrix $\Sigma$.  Is it true that,  
for any weights $0\leq w_j\leq 1$, $ j=1,\ldots, n, \,
\sum_{j=1}^n w_j=1$,
\begin{equation} \label{e:mod1}
\sum_{j=1}^n w_j \frac{X_j}{|Y_j|}\eid W,
\end{equation}
where $W$ is a standard Cauchy? As was the case with
\eqref{e:pillai.meng}, the claim \eqref{e:mod1} holds if the covariance matrix
$\Sigma$ is a diagonal matrix, and also when all the
correlations implied by the covariance matrix $\Sigma$ equal 1. 
However, we will show that \eqref{e:mod1} is not true 
for every covariance matrix $\Sigma$.
The argument uses a simple property of
Cauchy densities that follows immediately from \eqref{e:dens.cauchy}:
\begin{equation} \label{e:simple.c}
f(x;\sigma)=\frac{f(0;\sigma)}{\pi^2f^2(0;\sigma)x^2+1}, \ x\in\bbr.
\end{equation}

It is enough to show that \eqref{e:mod1} does not, in general, hold
when $n=2$. A non-singular covariance matrix $\Sigma$ with equal
variances can be scaled so that its inverse matrix is
\begin{equation} \label{e:sigma.inv}
\Sigma^{-1} = \left[\begin{matrix}
1 & \theta \\
\theta & 1
\end{matrix}\right] \ \ \text{for $-1<\theta<1$.}
\end{equation}
Let $R_i=X_i/|Y_i|, \, i=1,2$. An elementary calculation for the
joint density of $(R_1,R_2)$ gives, in the parametrization
\eqref{e:sigma.inv},
\begin{align*}
  &f_{R_1,R_2}(r_1,r_2)\\
  =&\int_{-\infty}^{\infty}\int_{-\infty}^{\infty}
                         \frac{(1-\theta^2)|y_1y_2|}{(2\pi)^2}\cdot \\ 
&\quad\quad\quad\quad\quad\exp\bigg\{-\frac{1}{2}(y_1^2r_1^2+2\theta|y_1y_2|r_1r_2+y_2^2r_2^2)
-\frac{1}{2}(y_1^2+2\theta  y_1y_2+y_2^2)\bigg\}dy_1dy_2\\
 =&2\int_{0}^{\infty}\int_{0}^{\infty}\frac{(1-\theta^2)xy_2^3}{(2\pi)^2} 
\exp\bigg\{-\frac{1}{2}\bigg(x^2y_2^2(r_1^2+1)+2\theta
    xy_2^2(r_1r_2+1)+y_2^2(r_2^2+1)\bigg)\bigg\}dxdy_2 \\
 +&2\int_{0}^{\infty}\int_{0}^{\infty}\frac{(1-\theta^2)xy_2^3}{(2\pi)^2} 
\exp\bigg\{-\frac{1}{2}\bigg(x^2y_2^2(r_1^2+1)+2\theta
   xy_2^2(r_1r_2-1)+y_2^2(r_2^2+1)\bigg)\bigg\}dxdy_2 \\
=&\int_{0}^{\infty}\frac{(1-\theta^2)x}{\pi^2\big(x^2(r_1^2+1)+2\theta
   x(r_1r_2+1)+(r_2^2+1)\big)^2}dx\\
  &+ 
\int_{0}^{\infty}\frac{(1-\theta^2)x}{\pi^2\big(x^2(r_1^2+1)+2\theta x(r_1r_2-1)+(r_2^2+1)\big)^2}dx.
\end{align*}
Define  the convex combination of $R_1$ and $R_2$ to be $V:=w_1R_1+w_2R_2$, where $w_1\in [0,1]$ and $w_2=1-w_1$.  Then the density $g_V$ of $V$ is given by
\begin{align} \notag
g_V(v) =&\int_{-\infty}^{\infty}\int_{0}^{\infty}\frac{(1-\theta^2)w_2^3x}
{\pi^2\big(w_2^2x^2(r^2+1)+2\theta w_2x(r(v-w_1r)+w_2)+((v-w_1r)^2+w_2^2)\big)^2}dxdr
  \\&+\int_{-\infty}^{\infty}\int_{0}^{\infty}\frac{(1-\theta^2)w_2^3x}
{\pi^2\big(w_2^2x^2(r^2+1)+2\theta w_2x(r(v-w_1r)-w_2)+((v-w_1r)^2+w_2^2)\big)^2}dxdr \label{e:dens1} \\
=&\int_{0}^{\infty}\frac{(1-\theta^2)(w_2^2x^2-2\theta w_1w_2x+w_1^2)x}{2\pi\big((w_2^2x^2-2\theta w_1w_2x+w_1^2)
(x^2+2\theta x+1)+(1-\theta^2)v^2x^2\big)^{3/2}}dx \notag \\
&+\int_{0}^{\infty}\frac{(1-\theta^2)(w_2^2x^2-2\theta w_1w_2x+w_1^2)x}{2\pi\big((w_2^2 x^2-2\theta w_1w_2x+w_1^2)
(x^2-2\theta x+1)+(1-\theta^2)v^2x^2\big)^{3/2}}dx. \notag 
\end{align}
For $w_1\in(0,1)$ fixed, we will prove that
\begin{equation} \label{e:lim1}
  \lim_{v\to\infty} v^2g_V(v)=1/\pi,
\end{equation}
for all $\theta\in(0,1)$.  
Assuming that \eqref{e:lim1} is true, if
$g_V$ were a Cauchy density for any $\theta\in(0,1)$ 
it would follow from \eqref{e:simple.c} that  
\begin{equation*}  
g_V(0)=1/\pi  \ \text{for any $\theta\in(0,1)$.}
\end{equation*}
We will show that, in fact,
\begin{equation} \label{e:g.zero1}
g_V(0)>1/\pi  
\end{equation}
for all $\theta\in(0,\epsilon)$ for some $\epsilon>0$.  
This precludes the possibility that $g_V$ is a Cauchy density for all such $\theta$.

To establish \eqref{e:g.zero1}, it suffices to show that  for any $w_1\in(0,1)$, 
the derivative of the function $g_V(0)$ with respect to $\theta$ is  positive  at  $\theta=0$. 
Indeed, by \eqref{e:dens1}
\begin{align*} 
g_V(0) =&\int_{0}^{\infty}\frac{(1-\theta^2)(w_2^2x^2-2\theta w_1w_2x+w_1^2)x}{2\pi\big(w_2^2x^2-2\theta w_1w_2x+w_1^2\big)^{3/2}
\bigl(x^2+2\theta x+1\big)^{3/2}}dx   \\
&+\int_{0}^{\infty}\frac{(1-\theta^2)(w_2^2x^2-2\theta w_1w_2x+w_1^2)x}{2\pi\big(w_2^2 x^2-2\theta w_1w_2x+w_1^2\bigr)^{3/2}
\bigl(x^2-2\theta x+1\big)^{3/2}}dx,
\end{align*}
It is elementary that one may differentiate with respect to $\theta$ inside the
integrals to obtain
$$
\frac{\partial g_V(0)}{\partial \theta} \bigg|_{\theta=0}
=\int_{0}^{\infty} \frac{w_1w_2x^2}
{\pi(w_1^2+w_2^2x^2)^{3/2}(x^2+1)^{3/2}} \, dx>0,
$$
as promised. We conclude that \eqref{e:g.zero1} holds  and, hence \eqref{e:mod1} cannot be valid for $\theta$ in a neighborhood of zero.

It remains to prove \eqref{e:lim1}. We use the last equality of \eqref{e:dens1} and write
in the obvious notation 
\begin{align*}
v^2g_V(v) =v^2\left[ \int_1^\infty +\int_1^\infty\right]
+ v^2\left[ \int_0^1 +\int_0^1\right]
  =I_1(v)+I_2(v).
\end{align*}
We have, as $v\to\infty$, 
\begin{align*}
I_1(v) \sim& (1-\theta^2)w_2^2\int_1^\infty \frac{x^3v^2}{\pi\bigl(
             w_2^2x^4 +(1-\theta^2)v^2x^2\bigr)^{3/2}}\,dx \\
  =&w_2 \int_{w_2/(v(1-\theta^2)^{1/2})}^\infty \frac{1}{\pi(x^2+1)^{3/2}}\,dx
     \to w_2 \int_0^\infty \frac{1}{\pi (x^2+1)^{3/2}}\,dx = w_2/\pi.
\end{align*}
Since $I_2$ can be reduced to $I_1$ by a variable change $x\to 1/x$
and switching $w_1$ and $w_2$, this means that $I_2(v)\to w_1/\pi$ as
$v\to\infty$. 
Because $w_1+w_2=1$, this proves
\eqref{e:lim1}.

In summary, it is not true that, for any choice of
$(w_1,w_2)$ with $0<w_1<1, \, w_1+w_2=1$, 
the sum on the left side of 
\eqref{e:mod1} has a Cauchy distribution for all $\theta$ in an
arbitrarily small neighborhood of the origin.

\section{Transformation 2: Randomly Stopped Brownian
  Motions} \label{sec:BM} 

Recall that an $\bbr^n$-valued Brownian motion $\bigl\{\bigl( X_1(t),\ldots,
X_n(t)\bigr), \, t\geq 0\bigr\}$ is a continuous centered Gaussian
process that starts at the origin at time zero and has stationary and 
independent increments. The law of such a process is completely
determined by the law of its values at time 1, $\bigl( X_1(1),\ldots,
X_n(1)\bigr)$. The covariance function of a Brownian motion is
\begin{equation} \label{e:BM.cov}
{\rm cov}(X_i(s),X_j(t))= \min(s,t)\, {\rm cov}(X_i(1),X_j(1)), \
s,t\geq 0, \ i,j=1,\ldots, n.
\end{equation}
A Brownian motion is self-similar in the sense that, for any $c>0$,
\begin{equation} \label{e:BM.selfs}
\bigl\{\bigl( X_1(ct),\ldots, X_n(ct)\bigr), \, t\geq 0\bigr\} \eid 
\bigl\{\bigl( c^{1/2}X_1(t),\ldots, c^{1/2}X_n(t)\bigr), \, t\geq 0\bigr\}.
\end{equation}
If $\bigl\{ X(t), \, t\geq 0\}$ is a one-dimensional Brownian motion, 
$Y$ is an independent copy of $X(1)$ and
$X(Y^{-2})$ is the value of the Brownian motion $X(t)$ evaluated
at the random time $t=Y^{-2}$,  then \eqref{e:BM.selfs} implies that
\begin{equation} \label{e:normal.cauchy.2}
X\bigl(Y^{-2}\bigr)\eid X(1)/|Y|\eid X(1)/Y\eid W
\end{equation}
by \eqref{e:normal.cauchy}, where $W$ is a standard Cauchy. 

Again following \cite{pillai:meng:2016}, 
let $\bigl\{\bigl( X_1(t),\ldots,
X_n(t)\bigr), \, t\geq 0\bigr\},\ n>1,$ be a Brownian motion with positive
variances 
and $(Y_1,\ldots, Y_n)$ a centered normal random
vector with the same law as $\bigl( X_1(1),\ldots,
X_n(1)\bigr)$ and independent of it. 
Let their common covariance matrix be $\Sigma$.
Is it true that for any weights 
$0\leq w_j\leq 1$, $ j=1,\ldots, n, \, \sum_{j=1}^n w_j=1$,
and any $\Sigma$,
\begin{equation} \label{e:mod2}
\sum_{j=1}^n w_j X_j\bigl(Y_j^{-2}\bigr)\eid W,
\end{equation}
where $W$ is a standard Cauchy? If $\Sigma$ 
is diagonal, or if all the
correlations implied by $\Sigma$ equal 1, 
then \eqref{e:mod2} holds. However, we will show that
\eqref{e:mod2} is not true  in general. The argument is similar to the
one used in Section \ref{sec:absolute.val} and relies on the property
\eqref{e:simple.c} of the Cauchy distribution. 
Again we let $n=2$ and parametrize the inverse of $\Sigma$ as in
\eqref{e:sigma.inv}.

Let $V=w_1X_1\bigl( Y_1^{-2}\bigr) + w_2X_2\bigl( Y_2^{-2}\bigr)$,
with $w_1\in [0,1]$ and $w_2=1-w_1$. Let $g_V$ be the density of $V$. If
$f_{Y_1,Y_2}$ is the joint density of $(Y_1,Y_2)$, then by
\eqref{e:BM.cov}, 
\begin{equation} \label{e:gV.2}
g_V(v) =\int_{-\infty}^\infty \int_{-\infty}^\infty h_{y_1,y_2}(v)
f_{Y_1,Y_2}(y_1,y_2)\, dy_1dy_2,
\end{equation}
where $h_{y_1,y_2}$ is the density of a centered normal random
variable with variance
$$
\frac{w_1^2}{y_1^2(1-\theta^2)} -
\frac{2w_1w_2\theta}{\max(y_1^2,y_2^2)(1-\theta^2)} +
\frac{w_2^2}{y_2^2(1-\theta^2)}.
$$
The part of the integral in \eqref{e:gV.2} over the range
$|y_1|<|y_2|$ can be written in the form
\begin{align*}
&\int_{-1}^{1}\int_{-\infty}^{\infty}\frac{(1-\theta^2)|x|y_2^2}{(2\pi)^{3/2}(w_1^2-2w_1w_2\theta x^2+w_2^2x^2)^{1/2}}\cdot \\ 
&\hskip 0.5in\exp\bigg\{-\frac{(1-\theta^2)x^2y_2^2l^2}{2(w_1^2-2w_1w_2\theta x^2+w_2^2x^2)}
-\frac{1}{2}(x^2y_2^2+2\theta x|y_2|y_2+y_2^2)\bigg\}dy_2dx\\
=&\int_{0}^{1}\frac{(1-\theta^2)(w_1^2-2w_1w_2\theta x^2+w_2^2x^2)x}
{2\pi\big((w_1^2-2w_1w_2\theta
   x^2+w_2^2x^2)(x^2+2\theta x+1)+(1-\theta^2)x^2v^2\big)^{3/2}}dx\\
&+\int_{0}^{1}\frac{(1-\theta^2)(w_1^2-2w_1w_2\theta x^2+w_2^2x^2)x}
{2\pi\big((w_1^2-2w_1w_2\theta
   x^2+w_2^2x^2)(x^2-2\theta
   x+1)+(1-\theta^2)x^2v^2\big)^{3/2}}dx.
\end{align*}
Computing in an identical manner the part of the integral in \eqref{e:gV.2} over the range
$|y_2|<|y_1|$ gives
\begin{align} \label{e:gV2.final}
g_V(v)=&\int_{0}^{1}\frac{(1-\theta^2)(w_1^2-2w_1w_2\theta
         x^2+w_2^2x^2)x} 
{2\pi\big((w_1^2-2w_1w_2\theta
   x^2+w_2^2x^2)(x^2+2\theta x+1)+(1-\theta^2)x^2v^2\big)^{3/2}}dx\\
&+\int_{0}^{1}\frac{(1-\theta^2)(w_1^2-2w_1w_2\theta x^2+w_2^2x^2)x}
{2\pi\big((w_1^2-2w_1w_2\theta
   x^2+w_2^2x^2)(x^2-2\theta
   x+1)+(1-\theta^2)x^2v^2\big)^{3/2}}dx \notag \\
&+\int_{0}^{1}\frac{(1-\theta^2)(w_2^2-2w_1w_1\theta x^2+w_1^2x^2)x}
{2\pi\big((w_2^2-2w_1w_2\theta
   x^2+w_1^2x^2)(x^2+2\theta
   x+1)+(1-\theta^2)x^2v^2\big)^{3/2}}dx\notag \\
&+\int_{0}^{1}\frac{(1-\theta^2)(w_2^2-2w_1w_2\theta x^2+w_1^2x^2)x}
{2\pi\big((w_2^2-2w_1w_2\theta
   x^2+w_1^2x^2)(x^2-2\theta
   x+1)+(1-\theta^2)x^2v^2\big)^{3/2}}dx. \notag 
\end{align}
We will prove that \eqref{e:lim1} still holds, independently of $\theta$ in \eqref{e:sigma.inv} 
and independently of the weights $(w_1,w_2)$.  
Assuming that this is the case, the
claim that $g_V$ cannot be a Cauchy density for all choices of $\theta$
and $(w_1,w_2)$ will, as in Section \ref{sec:absolute.val}, 
 follow once we show that, for any $0<w_1<1$,
the function $g_V(0)$ has a strictly positive derivative in $\theta$
at $\theta=0$. By \eqref{e:gV2.final},
\begin{align*} 
g_V(0)=&\int_{0}^{1}\frac{(1-\theta^2)(w_1^2-2w_1w_2\theta
         x^2+w_2^2x^2)x} 
{2\pi(w_1^2-2w_1w_2\theta
   x^2+w_2^2x^2)^{3/2}(x^2+2\theta x+1)^{3/2}}dx\\
&+\int_{0}^{1}\frac{(1-\theta^2)(w_1^2-2w_1w_2\theta x^2+w_2^2x^2)x}
{2\pi(w_1^2-2w_1w_2\theta
   x^2+w_2^2x^2)^{3/2}(x^2-2\theta
   x+1)^{3/2}}dx   \\
&+\int_{0}^{1}\frac{(1-\theta^2)(w_2^2-2w_1w_1\theta x^2+w_1^2x^2)x}
{2\pi(w_2^2-2w_1w_2\theta
   x^2+w_1^2x^2)^{3/2}(x^2+2\theta
   x+1)^{3/2}}dx \\
&+\int_{0}^{1}\frac{(1-\theta^2)(w_2^2-2w_1w_2\theta x^2+w_1^2x^2)x}
{2\pi(w_2^2-2w_1w_2\theta
   x^2+w_1^2x^2)^{3/2}(x^2-2\theta
   x+1)^{3/2}}dx. 
\end{align*}
It is elementary that one may differentiate with respect to $\theta$ inside the
integrals to obtain
$$
\frac{\partial g_V(0)}{\partial \theta} \bigg|_{\theta=0}
= \int_0^1 \frac{w_1w_2x^3}{\pi(w_1^2+w_2^2x^2)^{3/2}(x^2+1)^{3/2}}\,
dx
+ \int_0^1 \frac{w_1w_2x^3}{\pi(w_2^2+w_1^2x^2)^{3/2}(x^2+1)^{3/2}}\,
dx>0,
$$
as promised.

It remains to prove \eqref{e:lim1}. We use \eqref{e:gV2.final} and we
write in the obvious notation
$$
v^2g_V(v) = I_1(v)+I_2(v)+I_3(v)+I_4(v).
$$
We have, as $v\to\infty$,
\begin{align*}
I_1(v) \sim& (1-\theta^2)w_1^2\int_0^1 \frac{xv^2}{2\pi(w_1^2
  +(1-\theta^2)v^2x^2)^{3/2}}\,dx \\
  =& w_1 \int_0^{v(1-\theta^2)^{1/2}/w_1}\frac{x}{2\pi(x^2+1)^{3/2}}\,dx
     \to w_1 \int_0^\infty\frac{x}{2\pi(x^2+1)^{3/2}}\,dx=w_1/(2\pi).
\end{align*}
Similarly, as $v\to\infty$, 
$$
I_2(v)\to w_1/(2\pi), \ \ I_3(v)\to w_2/(2\pi), \ \ I_4(v)\to w_2/(2\pi),
$$
and \eqref{e:lim1} follows from the fact that
$w_1+w_2=1$. 

In summary, \eqref{e:mod2} does not hold for all values of the parameters. 
We prove the stronger contrary result that, for any choice of
$(w_1,w_2)$ with $0<w_1<1, \, w_1+w_2=1$,
the sum on the left side of \eqref{e:mod2} does not
have a Cauchy distribution for all $\theta$ in an
arbitrarily small neighborhood of the origin.

\section{Conclusions} \label{sec:conclusions}

Two natural transformations of multivariate normal vectors considered here
and the similar transformation \cite{pillai:meng:2016} considered
all result in a Cauchy distribution if the covariance
matrix of the normal ingredients is diagonal, or if all the
correlations implied by the covariance matrix equal 1. 
However, the transformation \eqref{e:pillai.meng} \cite{pillai:meng:2016} considered 
produces a Cauchy distribution regardless of the normal covariance matrix.
By contrast, the transformations \eqref{e:mod1} and \eqref{e:mod2}
do not always produce a Cauchy distribution. This shows that the
mysteries of the connections between the normal laws and the Cauchy
laws remain to be understood. The result of \cite{pillai:meng:2016}
exhibits a family of laws in $\bbr^n$ with standard Cauchy marginals
that share with a subclass of multivariate Cauchy laws
\eqref{e:mult.cauchy} the property that any convex linear combination
of the coordinates is again a standard Cauchy. Since the two
transformations we have considered lack this property, the
correlations between jointly normal random
variables are not always overwhelmed by the heaviness of the marginal
tails. 
Future work must clarify the extent of
the phenomenon described by \cite{drton:xiao:2016} and \cite{pillai:meng:2016}.


\end{document}